\documentclass[12pt]{article}
\usepackage{amsmath,amsthm,amsfonts}
\theoremstyle{plain}

\theoremstyle{plain}

\theoremstyle{remark}

\theoremstyle{remark}

\theoremstyle{plain}

\normalbaselines
\begin {document}

\def\ab{I\times{I}}
\def\bc{\{r_ir_j\}_{1\le i<j<\infty}}
\def\cd{\sum_{i=1}^\infty a_ix_i}
\def\de{L_{\infty}}
\def\ef{\sum_{1\le i<j<\infty}a_{i,j}r_ir_j}
\def\fg{\sum_{k=1}^{\infty}}
\def\gh{[x_n]_{n=1}^\infty}
\def\hi{\{x_n\}_{n=1}^\infty}
\def\ij{n\in{{\mathbb N}}}
\def\jk{\sigma_k}
\def\kl{\Delta_j^{(k)}}
\def\lm{\int_0^1}
\def\no{\sum_{1\le i<j\le n}}
\def\mn{{\bf r}_{i,j}}
\def\op{\bigg\|}
\def\pr{\bigg\|}
\def\rs{1\le i<j<\infty}
\def\ls{{\bar{\cal R}}}
\def\lv{2^k<i<j\le 2^{k+1}}
\def\lg{\varepsilon}
\def\lx{(u-3-2v)k/2}
\def\st{\varphi_{-1/2}}
\def\tv{\log_2^{1/2-\varepsilon}{2/t}}
\def\xy{\theta_{i,j}}
\def\nx{\varphi_{\varepsilon}}
\def\wg{{\bar T}_{\theta}}

\def\cl{\sum_{l=1}^k}



\title{RADEMACHER CHAOS IN SYMMETRIC SPACES, II}
\author{\sc S.\,V. Astashkin}
\maketitle

\begin{center}
{\it Department of Mathematics, Samara State University, \\
ul. Akad. Pavlova 1, 443011 Samara, RUSSIA\\
{\tt E-mail: astashkn@ssu.samara.ru}}
\end{center}


\begin{abstract}
In this paper we study some properties of the orthonormal
system $\bc$ where $r_k(t)$ are Rademacher functions on $[0,1],$
 $k=1,2,\ldots$ This system is usually called {\it Rademacher
 chaos of order~$2$}. It is shown that a specific ordering of the chaos leads
to a basic sequence (possibly non-unconditional) in a wide class of symmetric
 functional spaces on $[0,1]$. Necessary and sufficient conditions on the space
 are found for	the basic sequence $\bc$ to possess the unconditionality
 property.
\end{abstract}
\def\nl{\vskip 0.2cm}

\section{Introduction}

\nl This paper is a
 continuation of \cite{[1]} where we started the study of Radema\-cher chaos in
 functional symmetric spaces (s.s.) on the segment $[0,1]$. Let us first
recall some definitions and  notations from \cite {[1]}.
\par
As usual,
$$
r_k(t)=\,{\rm sign}\,\sin{2^{k-1}{\pi}t}\,\,(k=1,2,\ldots)
$$
denotes the system of Rademacher
 functions on $I:=[0,1]$. The set of all real-valued functions $x(t)$ that can
 be represented in the form
$$
x(t)=\;\sum_{\rs} a_{i,j}r_i(t)r_j(t)\;\;\;(t\in [0,1])
$$
is called a {\it chaos} of order $2$ with respect to the system
 $\{r_k(t)\}$ ({\it Rademacher chaos of order 2}\,).  The same name is used,
 with no ambiguity, for the orthonormal system of functions $\bc$. In the
 sequel, as in \cite {[1]},  $H$ denotes the closure of $\de$ in the Orlicz
 space $L_M$ where $M(t)=\,e^t-1$.
\medskip

 In \cite{[1]}, we proved  the
 following.
\medskip

{\bf Theorem~A.} {\it Let $X$ be a symmetric space. Then the following
statements are equivalent:
\smallskip

$1)$ The system $\bc$ in $X$ is equivalent to the canonical basis in~$l_2$;
\smallskip

$2)$ A continuous imbedding  $H\;\subset\;X$ takes place.}
\medskip

In this paper,
we shall consider questions related to the unconditionality of Rademacher
chaos.	Our  main result is:
\medskip
\par
{\it The statements $1)$
 and $2)$ in Theorem A are equivalent to the next one:
\smallskip

$3)$ The system $\bc$ is an unconditional basic sequence in s.s.~$X$.}
\medskip

Let us recall the  meaning of the central notions above.
\medskip

{\bf Definition.}
A sequence $\hi$ of elements in   Banach space $X$ is called a {\it
 basic} sequence if it is a basis in its closed linear span $\gh$.
\medskip

As is well-known  (see for example \cite[p.2]{[2]}\,), the latter is equivalent
to the following two conditions:

1) $x_n\,\ne\,0$ for all $\ij$;

2) The family of projectors
$$
P_m\left(\cd\right)\;=\;\sum_{i=1}^m a_ix_i\;\;\;(m=1,2,\ldots),
$$
defined on $\gh$, is uniformly bounded. That
is, a constant $K>0$ exists such that for all $m,n\in{\mathbb N}$, $m<n$, and
 $a_i\in{\mathbb R}$, the following inequality holds:
$$
\op\sum_{i=1}^m
 a_ix_i\pr\;\le\;K\,\op\sum_{i=1}^n a_ix_i\pr.\leqno{(1)}
$$

One of the most important properties of a basic sequence is its unconditionality.
\medskip

{\bf Definition.} A basic sequence $\hi$ in a
 Banach space $X$ is said to be {\it unconditional} if, for any rearrangement
 $\pi$ of ${\mathbb N}$, the sequence $\{x_{\pi(n)}\}_{n=1}^\infty$ is also a basic
sequence in  $X$.
\medskip

This is equivalent, in particular, to the uniform boundedness
of the family of operators
$$
M_\theta\bigg(\cd\bigg )\;=\;\sum_{i=1}^\infty
\theta_ia_ix_i\;\;\;(\theta_i=\,\pm 1)
$$
which are defined on $\gh$~\cite[p.18]{[2]}, and the last
 means that there is a	constant $K_0$	such that for each $\ij$ and any couple
 of sequences of signs $\{\theta_i\}$ and real numbers $\{a_i\}$,
$$
\op\sum_{i=1}^n \theta_ia_ix_i\pr\;\le\;K_0\,\op\sum_{i=1}^n
a_ix_i\pr.\leqno{(2)}
$$

Finally, note that for	sequences of real numbers $(a_{i,j})_{\rs}$ we	use
the common notation
$$
\|(a_{i,j})\|_2\;:=\;\bigg(\sum_{\rs}
a_{i,j}^2\bigg )^{1/2}. $$

\section{Rademacher chaos as a basic sequence}

\nl The system $\{r_k\}_{k=1}^\infty$ and Rademacher chaos $\bc$, both are
 special subsystems of	 Walsh system $\{w_n\}_{n=0}^\infty$. If the latter
 is considered with Paley indexing \cite[p.158]{[3]},
 then $w_0=r_1,\,w_{2^k}=r_{k+2},\;k=0,1,\ldots$ We shall enumerate Rademacher
 chaos in correspondence  to this indexing,
namely,
$$
\begin{array}{c} \varphi_1=\,r_1r_2=r_2,\,\varphi_2=\,r_1r_3=r_3,\,
\varphi_3=\,r_2r_3,\, \varphi_4=\,r_1r_4=r_4,\,\ldots,\\
\\
\varphi_{k(k-1)/2+1}=\,r_1r_{k+1}=r_{k+1},\,\ldots,\,\varphi_{k(k+1)/2}=\,
r_kr_ { k + 1 } ,\,\ldots
\end{array} \leqno{(3)}
$$

Before formulating our
first theorem let us recall the definition of a fundamental notion in the
 interpolation theory of operators (for more details, see
 \cite{[4]}).
\medskip

{\bf Definition.} A Banach space $X$ is said to be an {\it
 interpolation space} with respect to the Banach couple $(X_0,X_1)$ if $X_0\cap
 X_1\,\subset\, X\,\subset\,X_0+X_1$ and, in addition, if the boundedness of a
 linear operator $T$ in both $X_0$ and $X_1$ implies its boundedness in $X$ as
 well.
\medskip

{\bf Theorem 1.} {\it The Rademacher chaos $\bc$, ordered
 according to rule  $(3)$, is a basic sequence in every
interpolation with respect to the couple  $(L_1,\de)$ s.s. $X$ on $[0,1]$ .
\medskip

Proof.} In the sequel we shall use the
following property of Walsh systems  (see  \cite [p.45]{[5]}).
Introduce the  Fourier-Walsh partial sum operators
$$
S_px(t)\;:=\;\sum_{i=0}^p\lm
 x(s)w_i(s)\,ds\,\,w_i(t)\quad(p=0,1,\ldots)
$$
and denote $\jk:=\,S_{2^k}$.
Set
$$
\kl\;:=\;\left((j-1)2^{-k},j2^{-k}\right)
\quad(k=0,1,\ldots;\,j=1,2,\ldots,2^k).
$$
 Then
$$
\jk x(t)\;=\;2^k\,\int_{\kl} x(u)\,du
$$
for each
 $t\in\kl$. In other words,  the operator $\jk$ coincides with the averaging
 operator over the system of dyadic intervals $\{\kl\}_{j=1}^{2^k}$. It is easy
 to see that such an operator is bounded in $L_1$ and $\de$, and more
precisely, its norm is equal to $1$ in both spaces.  Thus $\jk$ is bounded in
$X$ as	well. Therefore, there exists a constant  $B=B(X)>0$ such that
$$
\|\jk
 x\|_X\;\le\;B\,\|x\|_X\leqno{(4)}
$$
for all $x\in X$.

For given  natural
 numbers $m<n$, and  real numbers $a_1, \ldots , a_n$,
 set
$$
y(t)=\;\sum_{i=1}^n a_i\varphi_i(t),\;\;z(t)=\;\sum_{i=1}^m
 a_i\varphi_i(t).
$$

In the simplest case, when  $\varphi_m=r_{k+2}$ for some
 $k=0,1,\ldots$, the orthonormality of Walsh system yields $z=\,\jk y$.
 Therefore, taking into account (4), we get $\|z\|\;\le\;B\,\|y\|$. Thus, for
 $x_i=\,\varphi_i$, inequality (1) holds with a constant $K=B$.

Consider now
 the general case: for some $0\le k\le l$, $2\le p<k+2$, $2\le
 q<l+2$,
$$
z(t)\;=\;\jk y(t)\,+\,\sum_{j=2}^p
 b_jr_j(t)\,r_{k+2}(t)
$$
and
$$
y(t)\;=\;\sigma_l y(t)\,+\,\sum_{j=2}^q
 c_jr_j(t)\,r_{l+2}(t).
$$
There are two possibilities.

{\small\sc Case 1.} $k=l,\,p\le q$.

Set
$$
f(t)\;=\;z(t)\,-\,\jk y(t)\;=\;\sum_{j=2}^p
 b_jr_j(t)\,r_{k+2}(t),$$
$$
g(t)\;=\;y(t)\,-\,z(t)\;=\;\sum_{j=p+1}^q
 c_jr_j(t)\,r_{k+2}(t).$$
It follows from the definition of Rademacher
 functions that the absolute values of $u(t)=\,f(t)+g(t)$ and $v(t)=\,f(t)-g(t)
$ are equimeasurable. Then the symmetry of $X$ yields $\|u\|_X\,=\,\|v\|_X$.
 Since $f=(u+v)/2$ we have
$$
\|f\|_X\;\le\;\|u\|_X.\leqno{(5)}
$$

Taking into account that
$$
u(t)\;=\;\sum_{j=2}^q c_jr_j(t)\,r_{k+2}(t)\;=\;y(t)\,-\,\jk
 y(t),
$$
the estimations (4) and (5) imply
$$
\|f\|_X\;\le\;\|y\|_X\,+\,\|\jk y\|_X\;\le\;(B+1)\|y\|_X,
$$
and consequently,
$$
\|z\|_X\;\le\;\|\jk
 y\|_X\,+\,\|f\|_X\;\le\;(2B+1)\|y\|_X.\leqno{(6)}
$$

{\small\sc Case 2.} $\;k<l$.

Now set
$$
g(t)\;=\;\sigma_{k+1}
y(t)\,-\,z(t)\;=\;
\sum_{j=p+1}^{k+1}d_jr_j(t)\,r_{k+2}(t)\,+\,d_{k+2}r_{k+3}(t).
$$
Define $f,\,u$ and $v$ as above. Inequality (5) holds in this case also.
 Since
$$
u\;=\;\sigma_{k+1} y\,-\,\jk y,
$$
inequality (4) implies
$$
\|f\|_X\;\le\;\|\sigma_{k+1}y\|_X\,+\,\|\jk y\|_X\;\le\;2B\|y\|_X.
$$
We have now
$$
\|z\|_X\;\le\;\|\jk
 y\|_X\,+\,\|f\|_X\;\le\;3B\|y\|_X.\leqno{(6')}
$$
The definitions of the
 functions $z$ and $y$, together with inequalities (6) and (6'), yield that
 relation (1) holds true for the  Rademacher chaos which is ordered according
to (3). The theorem is proved.
\medskip

{\bf Remark 1.} The requirement for
the space $X$ to be an interpolation space with respect to the couple
$(L_1,\de)$ is
 not very restrictive. The most important  s.s. (Orlicz, Lorentz, Marcinkiewicz
 spaces and others) possess this property \cite[p.142]{[4]}. In addition, it is
 seen from the proof of the theorem that the above condition may
 be replaced by a weaker one: the boundedness in $X$ of the averaging operators
 corresponding to the dyadic partitionings of the interval
 $[0,1]$.

\section{Rademacher chaos as unconditional basic sequence}

\nl We go now further  to the study of the unconditionality of Rademacher chaos
in s.s. We have already mentioned that the main result in this paper
 amplifies Theorem A proved in \cite{[1]} and formulated in Section~1.
\medskip

{\bf Theorem 2.} {\it Let $X$ be s.s. on $[0,1]$. Then the
 following assertions are equivalent:
\smallskip

$1)$ The system $\bc$ in $X$
 is equivalent to the canonical basis in the space $l_2$, that is, there is a
 constant $C>0$ that depends only on the space $X;$ such that  for all real
 numbers $a_{i,j}\;(\rs)$,
$$
C^{-1}\,\|(a_{i,j})\|_2\;\le\;\op\ef\pr_X\;\le\;C\,\|(a_{i,j})\|_2.
\leqno{(7)}
$$
\smallskip

$2)$ A continuous imbedding $H\;\subset\;X$ takes
 place;
\smallskip

$3)$ The system $\bc$ is an unconditional basic sequence in
 $X$.}
\medskip

{\bf Remark 2.} The implication $1)\Rightarrow 3)$ is evident
 and the equivalence $1)\Leftrightarrow 2)$ is proved in \cite{[1]}. Thus, it
 suffices to prove the implication $3)\Rightarrow 1)$.
\medskip

First, we
prove a weaker assertion. Let $G$ denote the closure of $\de$ in the Orlicz
 space $L_N$ corresponding to the function $N(t)=e^{t^2}-1$.
\medskip

{\bf
 Proposition 1.} {\it Let the s.s. $X$ on $[0,1]$ be such that $X\supset G$ and
 the system $\bc$ is an unconditional basic sequence in $X$. Then the assertion
 $1)$ in Theorem 2 holds true, that is, the system $\bc$ in $X$ is equivalent
to the canonical basis in $l_2$.}
\medskip

For the proof we need a lemma that
 concerns spaces with a mixed norm. Let us recall the definition (for details,
 see \cite[p.400]{[6]}\,).
\medskip

{\bf Definition.} Let $X$ and $Y$ be s.s.
 on $[0,1]$. The {\it space with a mixed norm} $X[Y]$ is the set of all
 measurable  functions $x(s,t)$ on the square $\ab $ satisfying the
 conditions:

1) $x(\cdot,t)\in Y$ for almost all $t\in I$;

2) $\varphi_x(t)\,:=\,\|x(\cdot,t)\|_Y\in X.$

Define
$$
\|x\|_{X[Y]}=\;\|\varphi_x\|_X.
$$

Let $A=A(u)$ be a
 $N$-function on $[0,\infty)$. This means that $A$ is continuous, convex, and
 satisfies
$$
\lim_{u\to +0}\frac{A(u)}{u}\;=\;\lim_{u\to
 +\infty}\frac{u}{A(u)}\;=\;0.
$$
As usual, denote by $L_A$ the Orlicz space of
 all functions $x=x(t)$ measurable on $[0,1]$ and  having a finite
 norm,
$$
\|x\|_{L_A}\;:=\;\inf\Big\{\lambda>0:\,
 \int_0^1 A\left(\frac{|x(t)|}{\lambda}\right)\,dt\,\le\,1\Big\}.
$$
Finally, let $A^*$ be the $N$-function conjugated to the $N$-function $A$, that
 is,
$$
A^*(u)\,:=\;\sup\{uv\,-\,A(v):\,v\ge 0\}.
$$

{\bf Lemma 1.} {\it The following imbeddings take place
$$
\de [L_A]\;\subset\;L_A(\ab),$$
$$
L_A(\ab) \;\subset\;L_1[L_A],$$
where $X(\ab)$ denotes a s.s. on the square $\ab$.
\medskip

Proof.} If $x=x(s,t)\in\de[L_A]$, then, according to the
 definition of the norm in an Orlicz space, for almost $t\in [0,1]$ we
 have
$$
\lm A\left(\frac{|x(s,t)|}{C}\right)\,ds\;\le\;1
$$
where
 $C=\,\|x\|_{\de[L_A]}$. After integrating this inequality and applying
Fubini's theorem we get
$$
\lm\lm
 A\left(\frac{|x(s,t)|}{C}\right)\,ds\,dt\;\le\;1.
$$
Therefore $x\in L_A(\ab)$
 and $\|x\|_{L_A(\ab)}\,\le\,\|x\|_{\de[L_A]}$. The first imbedding is
 proved.

For the proof of the second imbedding we pass on to the dual
 spaces. Recall that the dual space $X'$ to the s.s. $X$ consists of all
 measurable functions $y=y(t)$ for which
$$
\|y\|_{X'}\;:=\;\sup\left\{\lm
 x(t)y(t)\,dt:\;\|x\|_{X}\le 1\right\}\;<\,\infty.
$$
We have already proved that
$$
\de[L_{A^*}]\;\subset\;L_{A^*}(\ab).
$$
Therefore, for the dual spaces
 we
have
$$
\left(L_{A^*}(\ab)\right)'\;\subset\;\left(\de[L_{A^*}]\right)'.
$$
Finally, since $(L_A)'=L_{A^*}$, $(A^*)^*=A$ \cite[p.146, p.22]{[7]}
 and $(X[Y])'=X'[Y']$ \cite[Th.3.12]{[8]}, it follows
that
$$
L_A(\ab)\;\subset\;L_1[L_A].
$$
Besides,
$$
\|x\|_{L_1[L_A]}\;\le\;C\,\|x\|_{L_A(\ab)}
$$
for a certain $C>0$.
\medskip

{\it Proof of Proposition 1.}
 Let $\mn(u)$ $(\rs)$ denote Rademacher functions, arbitrarily ordered by the
 couples $(i,j)$.

Since $X\supset G$, by the lemma, for any given $\ij$ and
 real numbers $a_{i,j}$, we get
$$
\lm\,\op\no a_{i,j}\mn(u)r_ir_j\pr_X\,du\;=$$
$$
=\;\left\|\,\op\no a_{i,j}\mn(u)r_i(t)r_j(t)\pr_{X(t)}\right\|_{L_1(u) }\;\le$$
$$
\le~C_1\,\left\|\,\op\no a_{i,j}\mn(u)r_i(t)r_j(t)\pr_{G(u)}\right\|_{\de(t)}$$
(this follows from the
fact that $G$ is a subspace of $L_N$ and thus, for the functions in $\de(\ab)$,
 the norms in the spaces $L_1[G]$ and $L_1[L_N]$ coincide and so do the norms
in the spaces $\de[G]$ and $\de[L_N]$). By Khintchine's inequality for the
space $G$ (see, for example,  \cite{[9]}), this
yields
$$
\lm\,\op\no
a_{i,j}\mn(u)r_ir_j\pr_X\,du\;\le\;C_2\,\|(a_{i,j})\|_2.\leqno{(8)}
$$
On the other hand, as  is known
 \cite[Ch.4]{[10]},
$$
\op\ef\pr_{L_p}\;\asymp\;\|(a_{i,j})\|_2
$$
for any
$p\in [1,\infty)$. (This means that a constant $C>0$ exists depending only on
 $p$ such
that
$$
C^{-1}\,\|(a_{i,j})\|_2\;\le\;\op\ef\pr_{L_p}\;\le\;C\,\|(a_{i,j})\|_2).
$$
Therefore, by the imbedding $X\subset L_1$ which holds for each s.s. $X$
 on $[0,1]$ \linebreak \cite[p.124]{[4]}, we get
$$
\op\ef\pr_X\;\ge\;C_3\,\|(a_{i,j})\|_2.\leqno{(9)}
$$
Thus, the
inequality
$$
\lm\,\op\no
a_{i,j}\mn(u)r_ir_j\pr_X\,du\;\ge\;C_3\,\|(a_{i,j})\|_2,\leqno{(10)},
$$
which is opposite to (8), holds always true.

By the
 assumptions, with a constant depending only on the space $X$, we have
$$
\op
\no a_{i,j}r_ir_j\pr_X\;\asymp\;\lm\,\op\no a_{i,j}\mn(u)r_ir_j\pr_X\,du
$$
for
 each $\ij$ and all real numbers $a_{i,j}$. In this way, the proposition
follows from  relations (8) and (10).
\medskip

In \cite{[11]} (see also
 \cite{[12]}\,) the notion of RUC (random unconditional con\-ver\-gence)-system
 was introduced. We shall give here an equivalent definition.
\medskip

{\bf
 Definition.} Let $X$ be a Banach space and let $X^*$ be its dual space. The
 biorthogonal system $(x_n,x_n^*)$, $x_n\in X$, $x_n^*\in X^*$ $(n=1,2,\ldots)$
 is said to be a {\it RUC-system}, if there exists a constant $C>0$ such
 that
$$
\lm\,\op\sum_{i=1}^n r_i(s)x_i^*(x)x_i\pr_X\,ds\;\le\;C\|x\|_X
$$
for
 any $\ij$ and all $x\in\gh$ ($r_i(s)$ are Rademacher
 functions).
\medskip

Inequalities (8) and (9) yield the following.
\medskip

{\bf Corollary 1.} {\it If s.s. $X\supset G$, then
 Rademacher chaos $\bc$ together with the basic coefficients is a RUC-system in
 $X$.}
\medskip

{\bf Corollary 2.} {\it For each s.s. $X$ the following
 assertions are equivalent:

$1)\;X\supset G;$
\smallskip

$2)\;\lm\,\op\sum_{\rs}
 a_{i,j}\mn(u)r_ir_j\pr_X\,du\; \asymp\;\|(a_{i,j})\|_2$.
\medskip

Proof.} The
 implication $1)\Rightarrow 2)$ follows from inequalities (8) and (10).

 Suppose now that  2) takes place and let $a_{i,j}=0$ $(i\ne 1)$. From
 the definition of Rademacher functions and the assumptions we
 get
$$
\lm\,\op\sum_{j=2}^\infty a_{1,j}{\it r}_{1,j}(u)r_1r_j\pr_X
 du\; =\;\op\sum_{j=2}^\infty a_{1,j}r_j\pr_X\;
 \le\;C\,\bigg(\sum_{j=2}^\infty a_{1,j}^2\bigg )^{1/2}.
$$
Therefore  (see
 \cite[p.134]{[13]} or \cite{[14]}\,) $X\supset G$.
\medskip

{\bf
Corollary 3.} {\it Suppose that s.s. $X\supset G$. Then, for any
 set $\{a_{i,j}\}_{\rs}$ of real numbers, there exists a set of signs
 $\{\xy\}_{\rs}$, $\xy=\pm 1$,	such that
$$
\op\sum_{\rs}\xy
 a_{i,j}r_ir_j\pr_X\;\asymp\;\|(a_{i,j})\|_2.
$$

Proof.} Corollary 2 yields
$$
\inf\bigg\{\op\sum_{\rs}\xy
a_{i,j}r_ir_j\pr_X:\,\xy=\pm  1\bigg \}\;
$$
$$
\le\;\lm\,\op\sum_{\rs}
 a_{i,j}\mn(u)r_ir_j\pr_X\,du\;\le\;C\|(a_{i,j})\|_2.
$$
Now the assertion
 follows from the fact that the opposite inequality holds always (see
 (10)).
\medskip

In order to prove Theorem 2 we need some more auxiliary assertions.

Let $\{n_k\}_{k=1}^\infty$ be an increasing sequence of natural
 numbers. Set \linebreak $t_k=\,2^{-n_{k+1}}$,
 $m_k=\,(n_{k+1}-n_k)(n_{k+1}-n_k-1)/2$ and
$$
y_k(t)=\;\sum_{n_k<i<j\le
 n_{k+1}} r_i(t)r_j(t)\;\;(k=1,2,\ldots).
$$

{\bf Lemma 2.} {\it Let $c_k>0$
 $(k=1,2,\ldots)$,
$$
y(t)=\;\fg c_ky_k(t)\;\;(t\in [0,1]).
$$
If $y^*(t)$ is a
 non-increasing rearrangement of the function $|y(t)|$ {\rm\cite[p.83]{[4]} },
 then
$$
y^*(t_k)\;\ge\;\cl m_lc_l.\leqno{(11)}
$$

Proof.} By the definition
of Rademacher functions $y_l(t)=m_l$ provided that $0<t<2t_k$ and $1\le l\le
k$. Therefore
$$
y(t)\;=\;\cl m_lc_l\,+\,\sum_{l=k+1}^\infty c_ly_l(t)\;\;(0<t<2t_k).
$$
Besides, there exists a set $E\subset (0,2t_k)$ of
 Lebesgue measure $|E|=t_k$  such that
$$
\sum_{l=k+1}^\infty
 c_ly_l(t)\;\ge\;0 \quad  \hbox {for } ~~ t\in E.
$$
Applying the previous equality we get
$$
y(t)\;\ge\;\cl c_lm_l \quad \hbox { for } ~~  t\in E.
$$

 Inequality (11) follows now from the definition of the rearrangement and the
 fact that $|E|=t_k$.
\medskip

The next assertion makes  Theorem 8 from
 \cite{[1]}  more precise.  We shall use here the same notations as in \cite
 {[1]}.

If $X$ is a s.s. on $[0,1]$, then $\ls(X)$ denotes a subspace of $X$
 consisting of all functions of the form
$$
x(t)=\;\sum_{\rs}
 a_{i,j}r_i(t)r_j(t),\;\;(a_{i,j})_{\rs}\in l_2.
$$
 For any arrangement of
 signs (that is, for any sequence $\theta=\{\xy\}_{\rs}$, $\xy=\pm 1$), we
 define the operator
$$
\wg x(t)=\;\sum_{\rs} \xy a_{i,j}r_i(t)r_j(t)
$$ on the
 subspace $\ls(X)$,
\medskip

{\bf Proposition 2.} {\it There exists an arrangement of signs
 $\theta=\{\xy\}_{\rs}$ such that for each $\lg\in (0,1/2)$ one can find a
 function $x\in\ls(\de)$ satisfying
$$
(\wg x)^*(t)\;\ge\;b\,\tv
$$
with a
 constant $b>0$ independent of $t\in (0,1/16]$.
\medskip

Proof.} By Theorem 6
 in \cite{[1]} (see also Lemma 3 there), for each $k=1,2,\ldots$ one can find
 $\xy=\pm 1$ $(\lv)$ such that the functions
$$
z_k(t)=\;\sum_{\lv}\xy r_i(t)r_j(t)
$$
satisfy
$$
\|z_k\|_\infty\;\asymp\;2^{3k/2}.\leqno{(12)}
$$
Set $x_k(t)=\,2^{-(3+2\lg)k/2}z_k(t)$ and
$$
x(t)=\;\fg x_k(t)\;=\;\sum_{\rs} a_{i,j}r_i(t)r_j(t),
$$
where
 $a_{i,j}=\,2^{-(3+2\lg)k/2}\xy$, if  $\lv,\,k=1,2,\ldots$, and
 $a_{i,j}=0$, otherwise. It follows from (12) that
$$
\|x\|_\infty\;\le\;C\,\fg 2^{-\lg k}\;=\;C/(2^\lg-1),
$$
that is, $x\in \ls(\de)$.

Let the arrangement
 of signs $\theta$ consist of the values $\xy$ just determined for
 $\lv,\,k=1,2,\ldots$, and arbitrary $\xy$ for the other couples $(i,j),\,i<j$.
 We get
$$
y(t)=\;\wg x(t)\;=\;\fg 2^{-(3+2\lg)k/2}\,y_k(t)
$$
where
$$
y_k(t)=\;\sum_{\lv} r_i(t)r_j(t)\;\;(k=1,2,\ldots).
$$

Next we apply Lemma 2 to the  case when
 $n_k=\,2^k$ and $c_k=\,2^{-(3+2\lg)k/2}$. Then clearly $t_k=2^{-2^{k+1}}$ and
   $m_k\ge 2^{2k-2}$.  Therefore, for each $k=1,2,\ldots$, we have
$$
y^*(t_k)\;\ge\;\frac{1}{4}\sum_{i=1}^k
2^{2i}2^{-(3+2\lg)i/2}\;
\ge\;2^{(1/2-\lg)k-2}\;\ge\;C_1\,\log_2^{1/2-\lg}(2/t_k).
$$
Given an arbitrary $t\in (0,1/16]$, one can find a $k\in{\mathbb N}$ so
 that $t_{k+1}<t\le t_k$. Taking into account the previous inequality, we
get
$$
y^*(t)\;\ge\;y^*(t_k)\;\ge\;C_1\,\log_2^{1/2-\lg}(2/t_k)\;\ge$$
$$
\ge\;4^{\lg-1/2}C_1\,\log_2^{1/2-\lg}(2/t_{k+1})\;\ge\;b\,\tv.$$
The proposition is
 proved.
\medskip

Recall that Marcinkiewicz s.s. $M(\varphi)$ ($\varphi(t)\ge
 0$ is a concave increasing function on $[0,1]$ ) consists of all measurable
 functions $x(s)$ having a finite
norm
$$
||x||_{M(\varphi)}\,:=\,\sup\left\{\frac{1}{\varphi(t)}\,\int_o^t\,x^*(s)\,ds:
\,0<t\le 1\right\}.
$$

{\bf Corollary 4.} {\it Let $X$ be s.s. on
 $[0,1]$ such that, for any arrangement of signs $\theta=\{\xy\}_{\rs}$, the
 operator $\wg$ is bounded in $\ls(X)$.
 Then
$$
X\;\supset\;\bigcup_{\lg\in(0,1/2)}M(\nx)
$$
where $M(\nx)$ is
 Marcinkiewicz space determined by $\nx(t)=\,t\tv$.
\medskip

Proof.} Making
use of Proposition 2, we can find an arrangement of signs $\theta$ such that
for a given number $\lg\in(0,1/2)$ and some function $x\in\ls(\de)$
\linebreak ($x$ depends on $\lg$),
$$
(\wg x)^*(t)\;\ge\;b\,\tv\;\;(0<t\le 1/16).
$$

By the fact that
 $X\supset\de$ for any s.s. $X$ on $[0,1]$ \cite[p.124]{[4]} and taking into
 account the assumption and the symmetry of $X$, we get
$$
{\bar
 x}_{\lg}(t):=\;\tv\,\in\,X.\leqno{(13)}
$$
Now it follows from  the relation
 \cite[p.156]{[4]}
$$
\|y\|_{M(\nx)}\;\asymp\;\sup\left\{y^*(t)\,\tv:\,0<t\le
 1\right\}.
$$
 that ${\bar x}_{\lg}(t)$ has maximal rearrangement in the
 space $M(\nx)$. Besides,  by virtue of  (13),	$X\supset M(\nx)$. The
corollary is proved.
\medskip

{\bf Proposition 3.} {\it An arrangement of
signs $\theta=\{\xy\}_{\rs}$ exists such that for each $\lg\in (-1/2,1/2)$ and
 any $\delta\in (0,1/4-\lg/2)$ one can find a function $x\in\ls(M(\nx))$
 satisfying
$$
(\wg x)^*(t)\;\ge\;d\,\log_2^{1/2+\delta}(2/t)\leqno{(14)}
$$
where the constant
 $d>0$ does not depend on $t\in (0,1/16]$.
\medskip

Proof.} Let
 $\theta=\{\xy\}_{\rs}$, $z_k$ and $y_k$ $(k=1,2,\ldots)$ be defined in the
same way as in the proof of Proposition 2. It is well-known \cite{[14]} that
 Marcinkiewicz space $M(\st)\,=\,M(t\log_2(2/t))$ coincides with the
 Orlicz space $L_M$, $M(t)=e^t-1$. Therefore, by Theorem A from Section 1, we
 have
$$
\|z_k\|_{M(\st)}\;\asymp\;\bigg(\sum_{\lv}
 |\xy|^2\bigg )^{1/2}\;
\asymp\;2^k\;\;(k=1,2,\ldots).\leqno{(15)} $$

For any
 $u\in (0,1)$ the space $M(\nx)$ with $\lg=(1-2u)/2$ is a space of the type $u$
 with respect to the couple $(\de,M(\st))$, that is, a constant $C>0$ exists
 such
that
$$
\|x\|_{M(\nx)}\;\le\;C\,\|x\|_\infty^{1-u}\|x\|_{M(\st)}^u\leqno{(16)}
$$
for all $x\in\de$. In fact, as we have already mentioned in the proof
of Corollary~4,
$$
\|x\|_{M(\nx)}\;\le\;C'\,\sup\{x^*(t)\log_2^{-u}(2/t):\,0<t\le 1\}\;\le$$
$$
\le\;C'\,\left[\sup\{x^*(t):\,0<t\le1\}\right]^{1-u}\;
\left[\sup\{x^*(t)\,\log_2^{-1}(2/t):\,0<t\le 1\}\right]^u\;\le$$
$$
\le\;C\,\|x\|_\infty^{1-u}\,\|x\|_{M(\st)}^u.$$

From
(12), (15) and (16) we
get
$$
\|z_k\|_{M(\nx)}\;\le\;D\,2^{(3-u)k/2}\;\;(k=1,2,\ldots).\leqno{(17)}
$$

For a given $v>0$ (to be determined later), set
$$
x_k(t)=\;2^{\lx}\,z_k(t)\;,\;x(t)=\;\fg x_k(t).
$$
By virtue of  (17),
 $x\in\ls(M(\nx))$. Applying Lemma 2 to the function
$$
y(t)=\;\wg x(t)\;=\;\fg 2^{\lx}\,y_k(t)
$$
for $n_k=\,2^k$, $c_k=\,2^{\lx}$, we
get
$$
y^*(t)\;\ge\;\frac{1}{4}\,\sum_{i=1}^k
2^{2i}\,2^{(u-2v-3)i/2}\;\ge$$
$$
\ge\;C_1'\,2^{(1+u-2v)k/2}\;\ge\;C_1\,\log_2^{(1+u-2v)/2}(2/t_k)\;\;(k=1,2,
\ldots).$$
In the same way as in the proof of  Proposition 2 we conclude that
$$
y^*(t)\;\ge\;D_1\,\log_2^{(1+u-2v)/2}(2/t)\;\;(0<t\le 1/16).
$$

Since
  $\delta<1/4-\lg/2$ by assumption, we can choose $v$ so
 that $$0<v<1/4-\lg/2-\delta.$$ Therefore
$(u-2v)/2>\delta$ and inequality (14)
 holds with some $b>0$.
\medskip

{\bf Corollary 5.} {\it Let $X$ be s.s. on
 $[0,1]$ such that $X\supset M(\nx)$ for some $\lg\in(-1/2,1/2)$. If for any
 arrangement of signs $\theta=\{\xy\}_{\rs}$ the operator $\wg$ is bounded in
 $\ls(X)$,
then
$$
X\;\supset\;\bigcup_{0<\delta<1/4-\lg/2}\,M(\varphi_{-\delta}).
$$}

The proof is similar to that of  Corollary 4.
\medskip

Now we are ready to prove our  main Theorem 2.
\medskip

{\it Proof of Theorem 2.} As it has been
already mentioned in Remark 2, it suffices to verify the implication
 $3)\Rightarrow 1)$.

If the system $\bc$ is unconditional in s.s. $X$, then,
 for each arrangement of signs $\theta$, the operator $\wg$ is bounded in
 $\ls(X)$.  In particular, we get by Corollary 4 that
 $X\,\supset\,M(\varphi_{1/5})$. Therefore, by Corollary 5, it follows that
 $X\,\supset\,M(\varphi_{-1/10})$. Since $M(\varphi_{-1/10})\,\supset\,G$,
then all the more, $X\,\supset\,G$. Finally,   applying Proposition 1 we
 conclude that the system $\bc$ is equivalent to the canonical basis in
 $l_2$. The theorem is proved.
\medskip

{\bf Remark 3.} Assertions analogous
to Theorems 1 and 2 are valid for the multiple Rademacher system
 $\{r_i(s)r_j(t)\}_{i,j=1}^\infty$ considered on the square $\ab$, $I=[0,1]$,
as well. This follows from the equivalence of the symmetric norms for the
series with respect to the systems $\bc$ and $\{r_i(s)r_j(t)\}_{i,j=1}^\infty$
 (see \cite{[15]}, and also
 \cite{[1]}\,).


\end{document}